\documentclass[10pt, a4paper]{article}
\usepackage{amssymb}
\usepackage{amsmath}
\usepackage{amsthm}
\usepackage{latexsym}
\usepackage[dvips]{epsfig}
\usepackage{mathrsfs}
\usepackage{bm}
\usepackage{stmaryrd}
\usepackage{authblk}

\theoremstyle{plain}

\newtheorem{theorem}{Theorem}

\newtheorem{corollary}{Corollary}

\font\tenscr=rsfs10 scaled1100
\font\sevenscr=rsfs7 % scaled \magstep1
\font\fivescr=rsfs5 % scaled \magstep1
\skewchar\tenscr='177
\skewchar\sevenscr='177
\skewchar\fivescr='177
\newfam\scrfam
\textfont\scrfam=\tenscr
\scriptfont\scrfam=\sevenscr
\scriptscriptfont\scrfam=\fivescr

\newcommand{\tensor}[3]{_{#1\phantom{#2}#3}^{\phantom{#1}#2}}
\newcommand{\half}{\frac{1}{2}}

\newcounter{mnote}

\def\Pc{P}

\def\Yc{Y}
\def\nablac{\hat{\nabla}}
\def\Gammac{\hat{\Gamma}}
\def\bc{\hat{b}}
\def\Sc{S}

\def\Pp{\rho}

\def\Yp{y}
\def\nablap{\check{\nabla}}
\def\Gammap{\check{\Gamma}}
\def\bp{\check{b}}
\def\Sp{\Sigma}

\begin{document}
\title{\textbf{A note on the coincidence of \\the projective and conformal Weyl tensors}}

\author[]{{ Christian L\"ubbe} \footnote{E-mail address:{\tt c.luebbe@ucl.ac.uk}}}

\affil[]{Department of Mathematics, University College London, Gower
  Street, London WC1E 6BT, UK}

\maketitle

\begin{abstract}
This article examines the coincidence of the projective and conformal Weyl tensors associated to a given connection $\nabla $. The connection may be a general Weyl connection associated to a conformal class of metrics $[g]$. The main result for $n \ge 4$ is that the Weyl tensors coincide iff $\nabla $ is the Levi-Civita connection of an Einstein metric.
\end{abstract}

\section{Introduction}

In 1918 Hermann Weyl introduced, what is now known as Weyl geometries \cite{Weyl1918}. He observed that the Riemann curvature has a conformally invariant component $C\tensor{ij}{k}{l}$, which he referred to as the conformal curvature. In \cite{Weyl1921} Weyl discussed both conformal and projective geometries and showed that analogously the Riemann curvature has a projectively invariant component $W\tensor{ij}{k}{l}$, referred to as the projective curvature. The idea has been extend to parabolic geometries, (see e.g. \cite{BEG}, \cite{CSbook}) and in the modern literature the invariant curvature component is simply referred to as the Weyl tensor or the Weyl curvature, with the type of geometry typically implied by the context. In this article we will be dealing with $C\tensor{ij}{k}{l}$ and $W\tensor{ij}{k}{l}$ simultaneously and we will refer to them as the conformal and projective Weyl tensors respectively.

In \cite{Nur12} Nurowski investigated when a given projective class of connections $ [ \nabla ] $ on $M$ includes a Levi-Civita connection of some metric $g$ on $M$. An algorithm to check the metrisability of a chosen projective structure was given. 
In proposition 2.5 of \cite{Nur12} it was shown that the projective and conformal Weyl tensors coincide if and only if the Ricci tensor of the Levi-Civita connection satisfies
\begin{equation}
\label{Nurowski condition}
M\tensor{abcd}{ef}{} R_{ef}=0
\end{equation}
where 
$$
M\tensor{abcd}{ef}{} = 2 g_{a[c}\delta^e_{d]}\delta^f_{b} + 2 g_{a[d} g_{c]b}g^{ef} + 2(n-1) g_{b[d}\delta^f_{c]}\delta^e_{a}.
$$
Corollary 2.6 of \cite{Nur12} deduces that the projective and conformal Weyl tensors of an Einstein metric are equal. As a comment Nurowski raised the question whether there are non-Einstein metrics, which satisfy condition \eqref{Nurowski condition}.

This article proves that this is not the case. In particular, for a given connection $\nabla$ on an $n\ge 4$ dimensional manifold the projective and conformal Weyl tensors associated to $\nabla$ only agree if $\nabla$ is the Levi-Civita connection of an Einstein metric. 
The problem is addressed in more generality by allowing for general Weyl connections. This generalisation is of interest, due to the fact that neither the Ricci curvature of a general Weyl connection nor the Ricci curvature of a projective connection need be symmetric. Hence the possibility exists that the two Weyl tensors agree when using a general Weyl connection that is not a Levi-Civita connection for a metric in $[g]$.

\section{Projective and conformal connection changes}

We define the tensors
\begin{equation*}
\Sp^{kl}_{ij} = \delta^k_i \delta^l_j + \delta^l_i \delta^k_j ,
\quad\quad
\Sc^{kl}_{ij} = \delta^k_i \delta^l_j + \delta^l_i \delta^k_j - g_{ij} g^{kl}
\end{equation*}
Two connections $\nabla$ and $\nablap $ are projectively related if there exists a 1-form $\bp_i $ such that the connection coefficients are related by
\begin{equation*}
\Gammap\tensor{i}{k}{j} = \Gamma\tensor{i}{k}{j} + \Sp^{kl}_{ij} \bp_l
\end{equation*}
We denote the class of all connections projectively related to $\nabla$ by $[\nabla]_p $.

Suppose further that $\nabla$ is related to the conformal class $[g]$. By this we mean that there exists a 1-form $f_i$ such that
\begin{equation}
\label{Weyl connection metric condition}
\nabla_i g_{kl} = -2f_i g_{kl}
\end{equation}
This holds for $g_{ij}$ iff it holds for any representative in $[g]$. Connections that satisfy \eqref{Weyl connection metric condition} are referred to as general Weyl connections of $[g]$. Note that the Levi-Civita connection of any representative in $[g]$ satisfies \eqref{Weyl connection metric condition}. However $\nabla$ need not be the Levi-Civita connection for a metric in $[g]$. 

The connections $\nabla$ and $\nablac $ are conformally related if there exists a 1-form $\bc_i $ such that the connection coefficients are related by
\begin{equation*}
\Gammac\tensor{i}{k}{j} = \Gamma\tensor{i}{k}{j} + \Sc^{kl}_{ij} \bc_l
\end{equation*}
We denote the class of all connections conformally related to $\nabla$ by $[\nabla]_c $. Observe that all connections in $[\nabla]_c $ satisfy \eqref{Weyl connection metric condition}.

\section{Decomposition of the Riemann curvature}

Given a connection $\nabla$ the Riemann and Ricci tensors are defined as
\begin{equation*}
2 \nabla_{[i} \nabla_{j]} v^k = R\tensor{ij}{k}{l} v^l, \quad \quad R_{jl} = R\tensor{kj}{k}{l}
\end{equation*}
The projective and conformal Schouten tensor are related to the Ricci tensor of $\nabla$ by \cite{BEG}, \cite{Fri03}
\begin{eqnarray*}
\Pp_{ij} &=& \frac{1}{n-1} R_{(jl)} +  \frac{1}{n+1} R_{[jl]}\\
\Pc_{ij} &=& \frac{1}{n-2} R_{(jl)} +  \frac{1}{n} R_{[jl]} - \frac{R_{kl}g^{kl}}{2(n-2)(n-1)} g_{ij}
\end{eqnarray*}

\noindent The Schouten tensors can be used to decompose the Riemann curvature as follows
\begin{equation}
\label{Riemann decomposition}
R\tensor{ij}{k}{l}  = W\tensor{ij}{k}{l} + 2\Sp_{l[i}^{km} \Pp_{j]m} = C\tensor{ij}{k}{l} + 2\Sc_{l[i}^{km} \Pc_{j]m} ,
\end{equation}
where $W\tensor{ij}{k}{l}$ and $C\tensor{ij}{k}{l}$ are the projective and conformal Weyl tensors respectively. Moreover the once contracted Bianchi identity $\nabla_k R\tensor{ij}{k}{l} =0 $ implies \cite{BEG} that
\begin{eqnarray}\label{proj_Bianchi}
\nabla_k W\tensor{ij}{k}{l} &=& 2(n-2) \nabla_{[i} \Pp_{j]l} = (n-2)\Yp_{ijl}\\
\label{conf_Bianchi}
\nabla_k C\tensor{ij}{k}{l} &=& 2(n-3) \nabla_{[i} \Pc_{j]l}  = (n-3)\Yc_{ijl}.
\end{eqnarray}
The tensor $\Yp_{ijl}$ and $\Yc_{ijl}$ are known as the Cotton-York tensors.

Under a connection change $\check{\nabla} = \nabla + \bp$ respectively  $\hat{\nabla} = \nabla + \bc$ the Schouten tensors transform as
\begin{eqnarray*}
\Pp_{ij} - \check{\Pp}_{ij} &=& \nabla_i b_j + \half \Sp^{kl}_{ij} \bp_k \bp_l \\
\Pc_{ij} - \check{\Pc}_{ij} &=& \nabla_i b_j + \half \Sc^{kl}_{ij} \bc_k \bc_l 
\end{eqnarray*}
In both cases the Schouten tensors absorb all terms that arise in the Riemann tensor under connection changes. It follows that the projective Weyl tensor $W\tensor{ij}{k}{l} $ and the conformal Weyl tensor $C\tensor{ij}{k}{l} $ are invariants of the projective class $[\nabla]_p$ and the conformal class $[\nabla]_c$, respectively. The question we wish to address is for which manifolds these two invariants coincide.

\medskip
We note that for $n \le 2 \,$ $W\tensor{ij}{k}{l} =0$ and for $n \le 3 \,$ $C\tensor{ij}{k}{l} =0$. Therefore it follows trivially that:

\noindent \textit{In $n=2$ the Weyl tensors always agree. In $n=3$ they agree if and only if the manifold is projectively flat, i.e. the flat connection is contained in $[\nabla]_p $ }

Hence in the following we focus only on $n > 3$.

\section{Coincidence of the conformal and projective \\ Weyl tensors}
The Ricci tensor can be decomposed into its symmetric trace-free, skew and trace components with respect to the metric $g_{ij}$:
\begin{eqnarray}
\label{Riccidecomp}
R_{ij} &=& \Phi_{ij} + \varphi_{ij} + \frac{R}{n}g_{ij} %\\
\end{eqnarray}
Hence the Schouten tensors can be rewritten as
\begin{eqnarray}
\label{projectiveSchoutentoRicci}
\Pp_{ij} 
&=& \frac{1}{n-1} \Phi_{ij}  +  \frac{1}{n+1} \varphi_{ij} + \frac{R}{n(n-1)}g_{ij}\\
\label{conformalSchoutentoRicci}
\Pc_{ij} 
&=& \frac{1}{n-2} \Phi_{ij} + \frac{1}{n} \varphi_{ij} + \frac{R}{2n(n-1)}g_{ij} 
\end{eqnarray}
The condition $W\tensor{ij}{k}{l} = C\tensor{ij}{k}{l} $ is equivalent to
\begin{equation}
 2\Sp_{l[i}^{km} \Pp_{j]m} =  2\Sc_{l[i}^{km} \Pc_{j]m} 
\end{equation}

\noindent Substitutions of \eqref{projectiveSchoutentoRicci} and \eqref{conformalSchoutentoRicci}  give
\begin{eqnarray*}
2\Sp_{l[i}^{km} \Pp_{j]m} 
&=& \frac{2}{n-1}\delta_{[i}^k \Phi_{j]l}  + \frac{2 R}{n(n-1)}\delta_{[i}^k g_{j]l} + \frac{2}{n+1}\delta_{[i}^k \varphi_{j]l} - \frac{2}{n+1}  \delta_l^k \varphi_{ij}\\
2\Sc_{l[i}^{km} \Pc_{j]m} 
&=& \frac{2}{n-2}\delta_{[i}^k \Phi_{j]l} - \frac{2}{n-2} g_{l[i}\Phi_{j]m}g^{km} + \frac{2 R}{n(n-1)}\delta_{[i}^k g_{j]l} \nonumber \\
&& + \frac{2}{n}\delta_{[i}^k \varphi_{j]l}  - \frac{2}{n} g_{l[i}\varphi_{j]m}g^{km} - \frac{2}{n}  \delta_l^k \varphi_{ij}
\end{eqnarray*}
We observe that the scalar curvature terms are identical on both sides and hence only $\Phi_{ij}$ and $\varphi_{ij}$ are involved in our condition. The scalar curvature can take arbitrary values.

\noindent Taking the trace over $il$ and equating both sides.
\begin{eqnarray*}
2\Sp_{l[i}^{km} \Pp_{j]m} g^{il} 
&=& \frac{1}{n-1}\Phi\tensor{j}{k}{}  - \frac{R}{n} \delta_{j}^k + \frac{3}{n+1}\varphi\tensor{j}{k}{}  \\
2\Sp_{l[i}^{km} \Pc_{j]m} g^{il} 
&=& - \Phi\tensor{j}{k}{} - \frac{R}{n} \delta_{j}^k + \frac{4-n}{n}\varphi\tensor{j}{k}{} = - R\tensor{j}{k}{} + \frac{4}{n}\varphi\tensor{j}{k}{}
\end{eqnarray*}
Comparing irreducible components we find that we require
\begin{eqnarray}
\frac{n}{n-1}\Phi\tensor{j}{k}{} = 0 \quad \mathrm{and} \quad
\frac{n^2-4}{n(n+1)}\varphi\tensor{j}{k}{} = 0
\end{eqnarray}
Thus under our assumption of $n > 3$, both $\Phi_{ij}$ and $\varphi_{ij}$ must vanish. It follows that the Ricci tensor is pure trace and hence $g$ is an Einstein metric. Note that the Bianchi identities \eqref{proj_Bianchi}, \eqref{conf_Bianchi} imply that $R$ is constant.

The result can be formulated as follows
\begin{theorem}
Let $\nabla$ be a connection related to the conformal class $[g]$.
\begin{itemize}
\item In $n=2$ the Weyl tensors always vanish and hence agree. 
\item In $n=3$ the Weyl tensors  agree if and only if the manifold is projectively flat, i.e. the flat connection is contained in $[\nabla]_p $ 
\item In $n \ge 4$ the Weyl tensors agree if and only if the connection $\nabla$ is the Levi-Civita connection of the metric $g$ and the manifold is an Einstein manifold.
\end{itemize}
\end{theorem}

\begin{corollary}
If the projective and conformal Weyl tensor for $n\ge 4$ coincide then the Cotton-York tensors coincide as well. In fact they vanish identically.
\end{corollary}
The result follows immediately from the fact that the connection is the Levi-Civita connection of an Einstein metric. Hence the Schouten tensors are proportional to the metric and both Cotton-York tensors vanish.

\section{Conclusion}

It has been shown that the coincidence of the projective and conformal Weyl tensors is closely linked to the concept of Einstein metrics. For metric connections in $[\nabla]_c$ one could have deduced the main result directly from \eqref{Nurowski condition} by using the above decomposition of the Ricci tensor and using suitable traces of \eqref{Nurowski condition}. However, the set-up given here allowed for a direct generalisation to Weyl connections without requiring a more general form of \eqref{Nurowski condition}. Moreover it was felt that the set-up provided more clarity of role of the different types of curvatures involved.

\end{document}